\documentclass{article}

\usepackage{verbatim}
\usepackage{graphicx}
\usepackage{natbib}
\usepackage{caption}
\usepackage{algorithm}
\usepackage{algorithmic}
\usepackage{newfloat}
\usepackage{listings}
\usepackage{amsmath}
\usepackage{amssymb}
\usepackage{amsfonts}
\usepackage{amsopn}
\usepackage{epsfig}
\usepackage{microtype}
\usepackage{subfigure}
\usepackage{booktabs}
\usepackage{amsfonts}
\usepackage{nicefrac}
\usepackage{microtype}
\usepackage{xcolor}
\usepackage{stmaryrd}
\usepackage{lipsum}
\usepackage{rotating}
\usepackage{wasysym}
\usepackage{phonetic}
\usepackage{halloweenmath}
\usepackage{bm}
\usepackage{slashed}
\usepackage{amsmath,stackengine,scalerel}
\usepackage{bbm}
\usepackage{censor}
\usepackage[makeroom]{cancel}
\usepackage{mathtools}
\usepackage{amsthm}
\usepackage{tikz}
\usepackage{bigints}
\usepackage{comment}
\usepackage{accents}
\usepackage{trimclip}

\usepackage{url}

\newcommand{\tsp}[0]{{\rm T}}
\newcommand{\He}[0]{{H\!e}}

\newcommand{\mkron}[0]{{\otimes}}

\newcommand{\bg}[1]{{\bm{#1}}}

\title{A General Multiplication Theorem for Multivariate Hermite Polynomials}

\author{%
  Alistair Shilton \\
  Applied Artificial Intelligence Institute ${\rm (A}^2{\rm I}^2{\rm )}$ \\
  Deakin University, Geelong, Australia \\
  \texttt{alistair.shilton@deakin.edu.au} \\
}

\date{}

\begin{document}

\maketitle

\begin{abstract}
The multiplication theorem for univariate Hermite polynomials $H_k(\lambda x)$ 
is well-known.  In this paper we generalize this result to multivariate 
Hermite polynomials ${\rm H}_{\bf k}({\bg{\Lambda}}{\bf x};{\bg{\Sigma}})$, 
and use this result to derive a multiplication theorem for univariate 
polynomials applied to inner-products $H_k({\bg{\lambda}}^\tsp {\bf x})$.
\end{abstract}

\section{Introduction and Main Results}

Given symmetric, positive definite covariance ${\bg{\Sigma}} \in \mathbb{R}^{n 
\times n}$, the multivariate Hermite polynomials ${\rm H}_{\tt k} (\cdot; 
{\bg{\Sigma}}) : \mathbb{R}^n \to \mathbb{R}$, for multi-index ${\tt k} \in 
\mathbb{N}^n$, are \cite{Rah3}:
\[
 \begin{array}{rl}
  {\rm H}_{\tt k} \left( {\bf x}; {\bg{\Sigma}} \right)
  &\!\!\!\!=
  \left( -1 \right)^{|{\tt k}|}
  \exp \left( \frac{1}{2} {\bf x}^\tsp {\bg{\Sigma}}^{-1} {\bf x} \right)
  \left( \frac{\partial}{\partial {\bf x}} \right)^{\tt k}
  \exp \left( -\frac{1}{2} {\bf x}^\tsp {\bg{\Sigma}}^{-1} {\bf x} \right)
 \end{array}
\]
where $|{\tt k}| = \sum_i k_i$ and ${\bf a}^{\tt k} = \prod_i a_i^{k_i}$.  
Following convention, the multivariate probabilists' Hermite polynomials 
correspond to the special case $\He_{\tt k} (\cdot) = {\rm H}_{\tt k} (\cdot; 
{\bf I})$, and the physicists' Hermite polynomials $H_{\tt k} (\cdot) = {\rm 
H}_{\tt k} (\cdot; \frac{1}{2} {\bf I})$, with the univariate ($n=1$) cases 
denoted $H_k$ and $\He_k$ (with index $k \in \mathbb{N}$).  It is well known 
that univariate Hermite polynomials satisfy the multiplication property $\forall \lambda, x 
\in \mathbb{R}$: 
\begin{equation}
 \begin{array}{rl}
  \He_k \left( \lambda x \right)
  &\!\!\!\!=
  \mathop{\sum}\limits_{i=0}^{\lfloor\frac{k}{2}\rfloor}
  \frac{\left( \lambda^2-1 \right)^i}{i!2^i}
  \frac{k!}{(k-2i)!}
  \lambda^{k-2i}
  \He_{k-2i} \left( x \right) \\
  H_k \left( \lambda x \right)
  &\!\!\!\!=
  \mathop{\sum}\limits_{i=0}^{\lfloor\frac{k}{2}\rfloor}
  \frac{\left( \lambda^2-1 \right)^i}{i!}
  \frac{k!}{(k-2i)!}
  \lambda^{k-2i}
  H_{k-2i} \left( x \right) \\
 \end{array}
\label{eq:standard_mult}
\end{equation}
In this paper we show that, in the general case, for positive definite 
covariance matrices ${\bg{\Sigma}} \in \mathbb{R}^{n \times n}$, ${\bg{\Upsilon}} 
\in \mathbb{R}^{m \times m}$, for all ${\bg{\Lambda}} \in 
\mathbb{R}^{n \times m}, {\bf x} \in \mathbb{R}^m$, ${\tt k} \in \mathbb{N}^n$:
\begin{equation}
 \begin{array}{l}
  {\rm H}_{\tt k} \left( {{\bg{\Lambda}}}^\tsp {\bf x}; {\bg{\Sigma}} \right)
  =
  \mathop{\sum}\limits_{{{\tt q} \in \mathbb{N}^m : }\atop{|{\tt q}| \in \{|{\tt k}|,|{\tt k}|-2,\ldots\}}  }
  {\rm T}_{{\tt k},{\tt q}} \left( {{\bg{\Sigma}} {{\bg{\Lambda}}}^\tsp {\bg{\Upsilon}}^{-1}}; {\bg{\Sigma}}, {\bg{\Upsilon}} \right)
  {\rm H}_{{\tt q}} \left( {\bf x}; {\bg{\Upsilon}} \right)
 \end{array}
 \label{eq:mainresult}
\end{equation}
with coefficients:
\[
 \begin{array}{l}
  {\rm T}_{{\tt k},{\tt q}} \left( \tilde{\bg{\Lambda}}; {\bg{\Sigma}}, {\bg{\Upsilon}} \right)
  =
  \left.
  \frac{{\tt{k}!}}{2^i{\tt q}!i!}
  {\bf I}_{n}^{\mkron {\tt k}\tsp}
  \left(
  \tilde{\bg{\Lambda}}^{\mkron {\tt q}\;}
  \otimes
  \left( {\rm vec} \left( \tilde{\bg{\Lambda}} {\bg{\Upsilon}}^{-1} \tilde{\bg{\Lambda}}^\tsp - {\bg{\Sigma}}^{-1} \right) \right)^{\otimes i}
  \right)
  \right|_{i=\frac{|{\tt k}|-|{\tt q}|}{2}}
 \end{array}
\]
where ${\bf A}^{\otimes k} = {\bf A} \otimes {\bf A} \otimes \overset{k\;{\rm 
times}}{\ldots}{\otimes}{\bf A}$, ${\bf A}^{\mkron{\tt k}} = \bigotimes_i {\bf 
A}_{i:}^{\otimes k_i}$, $\otimes$ is the Kronecker product, and ${\rm vec}$ is 
the columnwise vectorization. This leads naturally to a series of results of 
intermediate generality, the most practical of which is $\forall {\bg{\lambda}}, 
{\bf x} \in \mathbb{R}^m$:
\begin{equation}
 \begin{array}{rl}
  \He_k \left( {{\bg{\lambda}}}^\tsp {\bf x} \right)
  &\!\!\!\!=
  \mathop{\sum}\limits_{{{\tt q} \in \mathbb{N}^m : }\atop{|{\tt q}| \in \{k,k-2,\ldots\}}}
  {T\!e}_{k,{\tt q}} \left( {{\bg{\lambda}}} \right)
  \prod_j \He_{q_j} \left( x_j \right)
  \\
  
  H_k \left( {{\bg{\lambda}}}^\tsp {\bf x} \right)
  &\!\!\!\!=
  \mathop{\sum}\limits_{{{\tt q} \in \mathbb{N}^m : }\atop{|{\tt q}| \in \{k,k-2,\ldots\}}}
  T_{k,{\tt q}} \left( {{\bg{\lambda}}} \right)
  \prod_j {H}_{q_j} \left( x_j \right)
 \end{array}
\label{eq:mainsubresult}
\end{equation}
where:
\[
 \begin{array}{rl}
  {T\!e}_{k,{\tt q}} \left( {{\bg{\lambda}}} \right)
  &\!\!\!\!=
  \left.
  \frac{{k!}}{2^i{\tt q}!i!}
  {\bg{\lambda}}^{{\tt q}}
  \left( \left\| {\bg{\lambda}} \right\|_2^2 - 1 \right)^i
  \right|_{i=\frac{k-|{\tt q}|}{2}} \\

   T_{k,{\tt q}} \left( {{\bg{\lambda}}} \right)
  &\!\!\!\!=
  \left.
  \frac{{k!}}{{\tt q}!i!}
  {\bg{\lambda}}^{{\tt q}}
  \left( \left\| {\bg{\lambda}} \right\|_2^2 - 1 \right)^i
  \right|_{i=\frac{k-|{\tt q}|}{2}}
\end{array}
\]
which may be helpful when applying Hermite transforms on 
vector projections.



The results presented here were derived in the course of other research but 
found to be superfluous for their original purpose. They are novel to the 
best of our knowledge, and we present them as a reference for others who 
recquire them.


\subsection{Notation} \label{sec:notation}

We use $\mathbb{N} = \{ 0,1,2,\ldots \}$.  
The number of elements in a finite set $\mathbb{S}$ is written $|\mathbb{S}|$. 
${\rm H}_{\tt k} (\cdot; {\bg{\Sigma}})$ are the multivariate Hermite polynomials.  
$\He_{\tt k}$ and $H_{\tt k}$ are the multivariate probabilists' and physicists' Hermite polynomials.  
$\He_k$ and $H_k$ are the univariate probabilists' and physicists' Hermite polynomials.  

{\bf Vectors and matrices: }
column vectors are denoted ${\bf a}, {\bf b}, \ldots$ and matrices ${\bf A}, {\bf B}, \ldots$. 
Vector and matrix elements are $a_i,b_j,\ldots$ and $A_{i,j}, B_{i,j}, \ldots$. 
Matrix rows and columns are ${\bf A}_{i:}$ and ${\bf A}_{:j}$.  
${\bf I}_n$ is the $n \times n$ identity matrix.  
${\bf A} \otimes {\bf B}$ is the Kronecker product. 
${\rm vec} ({\bf A})$ is the columnwise vectorization.
$[{\bf a}_i]_{:{\tt k}} = [{\bf a}_{k_1}; {\bf a}_{k_2}; \ldots]$ is a matrix constructed from column vectors ${\bf a}_i$. 

{\bf Multi-index notation:} 
multi-indices are denoted ${\tt k}, {\tt l}, \ldots$, 
with elements $k_i,l_j, \ldots$, 
and are assumed to be consistently (if arbitrarily) ordered.  
We use the notations 
$|{\tt k}| = \sum_i k_i$ (summation), 
${\bf x}^{\tt k} = \prod_i x_i^{k_i}$ (exponentiation), 
$(\frac{\partial}{\partial {\bf x}})^{\tt k} = \prod_i \frac{\partial^{k_i}}{\partial x_i^{k_i}}$ (derivatives), 
and ${\tt k}! = \prod_i k_i!$ (factorial).

{\bf Kronecker Powers:} 
${\bf A}^{\otimes k} = {\bf A} \otimes {\bf A} \otimes \overset{k\;{\rm times}\;}{\ldots}{\otimes}{\bf A}$ is the matrix Kronecker power and ${\bf a}^{\otimes k} = {\bf a} \otimes {\bf a} \otimes \overset{k\;{\rm times}\;}{\ldots}{\otimes}{\bf a}$ the vector Kronecker power, where we adopt the conventions ${\bf A}^{\otimes 0} = [1]$ and ${\bf a}^{\otimes 0} = [1]$; and 
${\bf A}^{\mkron {\tt k}} = \bigotimes_i {\bf A}_{i:}^{\otimes k_i}$ is the columnwise Kronecker power, where we note that 
this definition satisfies:
\begin{equation}
 \begin{array}{rl}
  ({\bf A}^\tsp {\bf b})^{\tt k} &\!\!\!\!= {\bf A}^{\mkron {\tt k}\tsp} {\bf b}^{\otimes |{\tt k}|}
 \end{array}
 \label{eq:kronprodexp}
\end{equation}

\section{Derivation of Main Result}

We begin by recalling the definition of the multivariate Hermite polynomials 
\cite[Definition 3]{Rah3}, given symmetric, positive definite ${\bg{\Sigma}} 
\in \mathbb{R}^{n \times n}$, indexed by ${\tt k} \in \mathbb{N}^n$:
\[
 \begin{array}{rl}
  {\rm H}_{\tt k} \left( {\bf x}; {\bg{\Sigma}} \right)
  &\!\!\!\!=
  \left( -1 \right)^{|{\tt k}|}
  \exp \left( \frac{1}{2} {\bf x}^\tsp {\bg{\Sigma}}^{-1} {\bf x} \right)
  \left( \frac{\partial}{\partial {\bf x}} \right)^{\tt k}
  \exp \left( -\frac{1}{2} {\bf x}^\tsp {\bg{\Sigma}}^{-1} {\bf x} \right)
  \;\;
  \forall {\bf x} \in \mathbb{R}^n
 \end{array}
\]
where we identify the special cases (respectively the probabilists' and 
physicists' multivariate Hermite polynomials):
\[
 \begin{array}{rl}
  \He_{\tt k} \left( {\bf x} \right)
  &\!\!\!\!=
  {\rm H}_{\tt k} \left( {\bf x}; {\bf I}_n \right) 
  =
  \left( -1 \right)^{|{\tt k}|}
  \exp \left( \frac{1}{2} {\bf x}^\tsp {\bf x} \right)
  \left( \frac{\partial}{\partial {\bf x}} \right)^{\tt k}
  \exp \left( -\frac{1}{2} {\bf x}^\tsp {\bf x} \right)
  =
  \prod_i
  \He_{k_i} \left( x_i \right) \\

  H_{\tt k} \left( {\bf x} \right)
  &\!\!\!\!=
  {\rm H}_{\tt k} \left( {\bf x}; \frac{1}{2} {\bf I}_n \right) 
  =
  \left( -1 \right)^{|{\tt k}|}
  \exp \left( {\bf x}^\tsp {\bf x} \right)
  \left( \frac{\partial}{\partial {\bf x}} \right)^{\tt k}
  \exp \left( -{\bf x}^\tsp {\bf x} \right)
  =
  \prod_i
  H_{k_i} \left( x_i \right)
 \end{array}
\]
where $\He_k$ and $H_k$ (for $k \in \mathbb{N}$) are the usual, univariate 
probabilists' and physicists' Hermite polynomials, respectively.

Noting that the multivariate Hermite polynomials satisfy 
\cite[Proposition 6]{Rah3}:
\begin{equation}
 \begin{array}{rl}
  \mathop{\sum}\limits_{{\tt k} \in \mathbb{N}^n}
  \frac{{\bf t}^{\tt k}}{{\tt k}!}
  {\rm H}_{\tt k} \left( {\bf x}; {\bg{\Sigma}} \right)
  &\!\!\!\!=
  \exp \left( {\bf t}^\tsp {\bg{\Sigma}}^{-1} {\bf x} - \frac{1}{2} {\bf t}^\tsp {\bg{\Sigma}}^{-1} {\bf t} \right)
 \end{array}
\label{eq:crux}
\end{equation}
we may proceed with our derivation.\footnote{We use as a template the proof of 
Theorem A.1 in \url{https://ncatlab.org/nlab/files/HermitePolynomialsAndHermiteFunctions.pdf}.}  
Throughout we assume ${\bg{\Sigma}} \in \mathbb{R}^{n \times n}$, 
${\bg{\Upsilon}} \in \mathbb{R}^{m \times m}$ are symmetric and 
positive-definite; and ${{\bg{\Lambda}}} \in \mathbb{R}^{m \times n}$, ${\bf x} 
\in \mathbb{R}^m$.

We begin by expanding the summation in (\ref{eq:crux}):
\[
 \begin{array}{l}
  \mathop{\sum}\limits_{{\tt k} \in \mathbb{N}^n}
  \frac{{\bf t}^{\tt k}}{{\tt k}!}
  {\rm H}_{\tt k} \left( {{\bg{\Lambda}}}^\tsp {\bf x}; {\bg{\Sigma}} \right)
  =
  \exp
  \left(
  {\bf t}^\tsp {\bg{\Sigma}}^{-1} {{\bg{\Lambda}}}^\tsp {\bf x}
  -
  \frac{1}{2} {\bf t}^\tsp {\bg{\Sigma}}^{-1} {\bf t}
  \right) \\


  =
  \exp
  \left(
  \tilde{\bf t}^\tsp {\bg{\Upsilon}}^{-1} {\bf x}
  -
  \frac{1}{2} \tilde{\bf t}^\tsp {\bg{\Upsilon}}^{-1} \tilde{\bf t}
  \right)
  \exp
  \left(
  \frac{1}{2} \tilde{\bf t}^\tsp {\bg{\Upsilon}}^{-1} \tilde{\bf t}
  -
  \frac{1}{2} {\bf t}^\tsp {\bg{\Sigma}}^{-1} {\bf t}
  \right) \\
 \end{array}
\]
where $\tilde{\bf t} = 
{\bg{\Upsilon}} {\bg{\Lambda}} {\bg{\Sigma}}^{-1} {\bf t}$. 
We observe that, again using (\ref{eq:crux}):
\[
 \begin{array}{l}
  \exp
  \left(
  \tilde{\bf t}^\tsp {\bg{\Upsilon}}^{-1} {\bf x}
  -
  \frac{1}{2} \tilde{\bf t}^\tsp {\bg{\Upsilon}}^{-1} \tilde{\bf t}
  \right)
  =
  \mathop{\sum}\limits_{{\tt q} \in \mathbb{N}^{m}}
  \frac{\left( 
  {\bg{\Upsilon}} {{\bg{\Lambda}}} {\bg{\Sigma}}^{-1} {\bf t} \right)^{{\tt q}}}{{\tt q}!}
  {\rm H}_{{\tt q}} \left( \sqrt{\frac{1}{s}\frac{m}{n}} {\bf x}; {\bg{\Upsilon}} \right) \\

  \exp
  \left(
  \frac{1}{2} \tilde{\bf t}^\tsp {\bg{\Upsilon}}^{-1} \tilde{\bf t}
  -
  \frac{1}{2} {\bf t}^\tsp {\bg{\Sigma}}^{-1} {\bf t}
  \right)
  =
  \exp
  \left(
  \frac{1}{2}
  {\bf t}^\tsp
  \left(
  {\bg{\Sigma}}^{-1} {{\bg{\Lambda}}}^\tsp {\bg{\Upsilon}} {\bg{\Upsilon}}^{-1} {\bg{\Upsilon}} {{\bg{\Lambda}}} {\bg{\Sigma}}^{-1}
  -
  {\bg{\Sigma}}^{-1}
  \right)
  {\bf t}
  \right) \\
 \end{array}
\]
which we rewrite in terms of (columnwise) Kronecker powers, using (\ref{eq:kronprodexp}):
\[
 \begin{array}{l}
  \exp
  \left(
  \tilde{\bf t}^\tsp {\bg{\Upsilon}}^{-1} {\bf x}
  -
  \frac{1}{2} \tilde{\bf t}^\tsp {\bg{\Upsilon}}^{-1} \tilde{\bf t}
  \right)
  =
  \mathop{\sum}\limits_{{\tt q} \in \mathbb{N}^{m}}
  \frac{\left( {\bg{\Sigma}}^{-1} {{\bg{\Lambda}}}^\tsp {\bg{\Upsilon}} \right)^{\mkron {\tt q}\tsp} {\bf t}^{\otimes |{\tt q}|}} {{\tt q}!}
  {\rm H}_{{\tt q}} \left( 
  {\bf x}; {\bg{\Upsilon}} \right) \\

  \exp
  \left(
  \frac{1}{2} \tilde{\bf t}^\tsp {\bg{\Upsilon}}^{-1} \tilde{\bf t}
  -
  \frac{1}{2} {\bf t}^\tsp {\bg{\Sigma}}^{-1} {\bf t}
  \right)
  =
  \mathop{\sum}\limits_{i \in \mathbb{N}}
  \frac{  \left(
  {\rm vec}
  \left(
  {\bg{\Sigma}}^{-1} {{\bg{\Lambda}}}^\tsp {\bg{\Upsilon}}
  {\bg{\Upsilon}}^{-1}
  {\bg{\Upsilon}} {{\bg{\Lambda}}} {\bg{\Sigma}}^{-1}
  -
  {\bg{\Sigma}}^{-1}
  \right)^{\otimes i}
  \right)^{\tsp}
  {\bf t}^{\otimes 2i}}{2^i i!}
 \\
 \end{array}
\]
and subsequently:
\[
 \begin{array}{l}
  \mathop{\sum}\limits_{{\tt k} \in \mathbb{N}^n}
  \frac{{\bf t}^{\tt k}}{{\tt k}!}
  {\rm H}_{\tt k} \left( {{\bg{\Lambda}}}^\tsp {\bf x}; {\bg{\Sigma}} \right)
  =
  \mathop{\sum}\limits_{{\tt q} \in \mathbb{N}^{m}}
  \mathop{\sum}\limits_{i\in\mathbb{N}}
  {\rm H}_{{\tt q}} \left( 
  {\bf x}; {\bg{\Upsilon}} \right)
  \frac{1}{2^i{\tt q}!i!}
  \ldots \\ \ldots
  \left(
  \left( 
  {\bg{\Sigma}}^{-1} {{\bg{\Lambda}}}^\tsp {\bg{\Upsilon}}
  \right)^{\mkron {\tt q}}
  \otimes
  \left( {\rm vec} \left( 
  {\bg{\Sigma}}^{-1} {{\bg{\Lambda}}}^\tsp {\bg{\Upsilon}}
  {\bg{\Upsilon}}^{-1}
  {\bg{\Upsilon}} {{\bg{\Lambda}}} {\bg{\Sigma}}^{-1}
  - {\bg{\Sigma}}^{-1}
  \right)
  \right)^{\otimes i}
  \right)^\tsp
  {\bf t}^{\otimes |{\tt q}|+2i} \\
 \end{array}
\]
From this, Recalling (\ref{eq:kronprodexp}), we have that $\forall {\bf{t}} \in \mathbb{R}^n$
\[
 \begin{array}{l}
  \mathop{\sum}\limits_{k \in \mathbb{N}}
  {\bf t}^{\otimes k\tsp}
  \mathop{\sum}\limits_{{{\tt k} \in \mathbb{N}^n :}\atop{ |{\tt k}| = k}}
  {\bf I}_{n}^{\mkron {\tt k}}
  \frac{1}{{\tt k}!}
  {\rm H}_{\tt k} \left( {{\bg{\Lambda}}}^\tsp {\bf x}; {\bg{\Sigma}} \right)
  =
  \mathop{\sum}\limits_{k \in \mathbb{N}}
  {\bf t}^{\otimes k\tsp}
  \mathop{\sum}\limits_{{{\tt q} \in \mathbb{N}^{m}, i\in\mathbb{N} :}\atop{|{\tt q}|+2i = k}}
  \frac{1}{2^i{\tt q}!i!}
  \ldots \\ \hfill 
  \left(
  \left(  
  {\bg{\Sigma}}^{-1} {{\bg{\Lambda}}}^\tsp {\bg{\Upsilon}}
  \right)^{\mkron {\tt q}}
  \otimes
  \left( {\rm vec} \left(
  {\bg{\Sigma}}^{-1} {{\bg{\Lambda}}}^\tsp {\bg{\Upsilon}}
  {\bg{\Upsilon}}^{-1}
  {\bg{\Upsilon}} {{\bg{\Lambda}}} {\bg{\Sigma}}^{-1}
  - {\bg{\Sigma}}^{-1} \right) \right)^{\otimes i}
  \right)
  {\rm H}_{{\tt q}} \left( 
  {\bf x}; {\bg{\Upsilon}} \right) \\
 \end{array}
\label{eq:stepaaa}
\]
Equating like-terms in the power series, this is equivalent to, $\forall k \in \mathbb{N}$:
\[
 \begin{array}{l}
  \left[ \begin{array}{c}
  {\bf I}_{n}^{\mkron {\tt k}}
  \end{array} \right]_{:{{\tt k} \in \mathbb{N}^n, |{\tt k}|=k}}
  \left[ \begin{array}{c}
  \frac{1}{{\tt k}!}
  {\rm H}_{\tt k} \left( {{\bg{\Lambda}}}^\tsp {\bf x}; {\bg{\Sigma}} \right)
  \end{array} \right]_{{{\tt k} \in \mathbb{N}^n, |{\tt k}|=k}}
  =
  \ldots \\ \hfill 
  \left[ \begin{array}{c}
  \frac{1}{2^i{\tt q}!i!}
  \left(
  \left(  
  {\bg{\Sigma}}^{-1} {{\bg{\Lambda}}}^\tsp {\bg{\Upsilon}}%
  \right)^{\mkron {\tt q}}
  \otimes
  \left( {\rm vec} \left( 
  {\bg{\Sigma}}^{-1} {{\bg{\Lambda}}}^\tsp {\bg{\Upsilon}}
  {\bg{\Upsilon}}^{-1}
  {\bg{\Upsilon}} {{\bg{\Lambda}}} {\bg{\Sigma}}^{-1}
  - {\bg{\Sigma}}^{-1} \right) \right)^{\otimes i}
  \right)
  \end{array} \right]_{:{{\tt q} \in \mathbb{N}^{m}, i\in\mathbb{N} :}\atop{|{\tt q}|+2i = k}}
  \\ \hfill\ldots
  \left[ \begin{array}{c}
  {\rm H}_{{\tt q}} \left( 
  {\bf x}; {\bg{\Upsilon}} \right)
  \end{array} \right]_{{{\tt q} \in \mathbb{N}^{m}, i\in\mathbb{N} :}\atop{|{\tt q}|+2i = k}}
  \;
 \end{array}
\]
Noting that $[ {\bf I}_{n}^{\mkron {\tt k}} ]_{:{{\tt k} \in \mathbb{N}^n, |{\tt k}|=k}}$ has 
$p = |\{{{\tt k} \in \mathbb{N}^n, |{\tt k}|=k}\}|$ orthonormal columns of dimension $n^k \geq p$, it follows that 
$[ {\bf I}_{n}^{\mkron {\tt k}} ]_{:{{\tt k} \in \mathbb{N}^n, |{\tt k}|=k}}^\tsp 
[ {\bf I}_{n}^{\mkron {\tt k}} ]_{:{{\tt k} \in \mathbb{N}^n, |{\tt k}|=k}} = {\bf I}_p$, so 
$\forall k \in \mathbb{N}$:
\[
 \begin{array}{l}
  \left[ \begin{array}{c}
  \frac{1}{{\tt k}!}
  {\rm H}_{\tt k} \left( {{\bg{\Lambda}}}^\tsp {\bf x}; {\bg{\Sigma}} \right)
  \end{array} \right]_{{{\tt k} \in \mathbb{N}^n, |{\tt k}|=k}}
  =
  \left[ \begin{array}{c}
  {\bf I}_{n}^{\mkron {\tt k}}
  \end{array} \right]_{:{{\tt k} \in \mathbb{N}^n, |{\tt k}|=k}}^\tsp
  \ldots \\ \hfill 
  \left[ \begin{array}{c}
  \frac{1}{2^i{\tt q}!i!}
  \left(
  \left(  
  {\bg{\Sigma}}^{-1} {{\bg{\Lambda}}}^\tsp {\bg{\Upsilon}}
  \right)^{\mkron {\tt q}}
  \otimes
  \left( {\rm vec} \left(
  {\bg{\Sigma}}^{-1} {{\bg{\Lambda}}}^\tsp {\bg{\Upsilon}}
  {\bg{\Upsilon}}^{-1}
  {\bg{\Upsilon}} {{\bg{\Lambda}}} {\bg{\Sigma}}^{-1}
  - {\bg{\Sigma}}^{-1} \right) \right)^{\otimes i}
  \right)
  \end{array} \right]_{:{{\tt q} \in \mathbb{N}^{m}, i\in\mathbb{N} :}\atop{|{\tt q}|+2i = k}}
  \\ \hfill\ldots
  \left[ \begin{array}{c}
  {\rm H}_{{\tt q}} \left( 
  {\bf x}; {\bg{\Upsilon}} \right)
  \end{array} \right]_{{{\tt q} \in \mathbb{N}^{m}, i\in\mathbb{N} :}\atop{|{\tt q}|+2i = k}}
  \;
 \end{array}
\]
which we can immediately rewrite componentwise $\forall {\tt k} \in \mathbb{N}^n$ as:
\[
 \begin{array}{l}
  {\rm H}_{\tt k} \left( {{\bg{\Lambda}}}^\tsp {\bf x}; {\bg{\Sigma}} \right)
  =
  \mathop{\sum}\limits_{{{\tt q} \in \mathbb{N}^m : }\atop{|{\tt q}| \in \{|{\tt k}|,|{\tt k}|-2,\ldots\}}}
  {\rm T}_{{\tt k},{\tt q}} \left( 
  {\bg{\Sigma}} {{\bg{\Lambda}}}^\tsp {\bg{\Upsilon}}^{-1}
  ; {\bg{\Sigma}}, {\bg{\Upsilon}} \right)
  {\rm H}_{{\tt q}} \left( 
  {\bf x}; {\bg{\Upsilon}} \right)
 \end{array}
\]
where:
\[
 \begin{array}{l}
  {\rm T}_{{\tt k},{\tt q}} \left( \tilde{\bg{\Lambda}}; {\bg{\Sigma}}, {\bg{\Upsilon}} \right)
  =
  \left.
  \frac{{\tt{k}!}}{2^i{\tt q}!i!}
  {\bf I}_{n}^{\mkron {\tt k}\tsp}
  \left(
  \tilde{\bg{\Lambda}}^{\mkron {\tt q}}
  \otimes
  \left( {\rm vec} \left( \tilde{\bg{\Lambda}} {\bg{\Upsilon}}^{-1} \tilde{\bg{\Lambda}}^\tsp - {\bg{\Sigma}}^{-1} \right) \right)^{\otimes i}
  \right)
  \right|_{i=\frac{|{\tt k}|-|{\tt q}|}{2}}
 \end{array}
\]
which is the result given (\ref{eq:mainresult}). Note that the 
standard multiplication property (\ref{eq:standard_mult}) for the 
univariate Hermite polynomials follows in the special case 
$n = m = 1$ with appropriate covariance matrices and subsequent 
simplification.


\subsection{Special Cases}

The main result (\ref{eq:mainresult}) 
is rather general, so it may be helpful to consider special cases of particular 
interest.  Suppose first that ${\bg{\Sigma}} = \sigma^2 {\bf I}_n$, ${\bg{\Upsilon}} 
= \sigma^2 {\bf I}_m$ (where $\sigma \in \mathbb{R}_+$). In this case 
(\ref{eq:mainresult}) becomes:
\[
 \begin{array}{l}
  {\rm H}_{\tt k} \left( {{\bg{\Lambda}}}^\tsp {\bf x}; \sigma^2 \right)
  =
  \mathop{\sum}\limits_{{{\tt q} \in \mathbb{N}^m : }\atop{|{\tt q}| \in \{|{\tt k}|,|{\tt k}|-2,\ldots\}}}
  {\rm T}_{{\tt k},{\tt q}} \left( {{\bg{\Lambda}}}; \sigma^2 \right)
  {\rm H}_{{\tt q}} \left( {\bf x}; \sigma^2 \right)
 \end{array}
\]
where we simplify the coefficients to:
\[
 \begin{array}{l}
  {\rm T}_{{\tt k},{\tt q}} \left( {{\bg{\Lambda}}} ; \sigma^2 \right)
  =
  \left.
  \frac{{\tt{k}!}}{(2\sigma^2)^i{\tt q}!i!}
  {\bf I}_{n}^{\mkron {\tt k}\tsp}
  \left(
  {\bg{\Lambda}}^{\mkron {\tt q}}
  \otimes
  \left( {\rm vec} \left( {\bg{\Lambda}} {\bg{\Lambda}}^\tsp - {\bf I}_{m} \right) \right)^{\otimes i}
  \right)
  \right|_{i=\frac{|{\tt k}|-|{\tt q}|}{2}}
 \end{array}
\]
In the case $n=1$, replacing ${\bg{\Lambda}}$ with ${\bg{\lambda}}$, 
this further simplifies to:
\[
 \begin{array}{l}
  {\rm H}_k \left( {{\bg{\lambda}}}^\tsp {\bf x}; \sigma^2 \right)
  =
  \mathop{\sum}\limits_{{{\tt q} \in \mathbb{N}^m : }\atop{|{\tt q}| \in \{k,k-2,\ldots\}}}
  {\rm T}_{k,{\tt q}} \left( {{\bg{\lambda}}}; \sigma^2 \right)
  {\rm H}_{{\tt q}} \left( {\bf x}; \sigma^2 \right)
 \end{array}
\]
where:
\[
 \begin{array}{l}
  {\rm T}_{k,{\tt q}} \left( {{\bg{\lambda}}}; \sigma^2 \right)
  =
  \left.
  \frac{{k!}}{(2\sigma^2)^i{\tt q}!i!}
  {\bg{\lambda}}^{{\tt q}}
  \left( \left\| {\bg{\lambda}} \right\|_2^2 - 1 \right)^i
  \right|_{i=\frac{k-|{\tt q}|}{2}}
 \end{array}
\]
Finally, for the special cases of probabilists' ($\sigma^2=1)$ and physicists' ($\sigma^2=\frac{1}{2}$) Hermite polynomials, 
we obtain:
\[
 \begin{array}{l}
  \He_k \left( {{\bg{\lambda}}}^\tsp {\bf x} \right)
  =
  \mathop{\sum}\limits_{{{\tt q} \in \mathbb{N}^m : }\atop{|{\tt q}| \in \{k,k-2,\ldots\}}}
  {T\!e}_{k,{\tt q}} \left( {{\bg{\lambda}}} \right)
  \He_{{\tt q}} \left( {\bf x} \right)
  \\
  
  H_k \left( {{\bg{\lambda}}}^\tsp {\bf x} \right)
  =
  \mathop{\sum}\limits_{{{\tt q} \in \mathbb{N}^m : }\atop{|{\tt q}| \in \{k,k-2,\ldots\}}}
  T_{k,{\tt q}} \left( {{\bg{\lambda}}} \right)
  {H}_{{\tt q}} \left( {\bf x} \right)
 \end{array}
\]
or, equivalently:
\[
 \begin{array}{l}
  \He_k \left( {{\bg{\lambda}}}^\tsp {\bf x} \right)
  =
  \mathop{\sum}\limits_{{{\tt q} \in \mathbb{N}^m : }\atop{|{\tt q}| \in \{k,k-2,\ldots\}}}
  {T\!e}_{k,{\tt q}} \left( {{\bg{\lambda}}} \right)
  \prod_j \He_{q_j} \left( x_j \right)
  \\
  
  H_k \left( {{\bg{\lambda}}}^\tsp {\bf x} \right)
  =
  \mathop{\sum}\limits_{{{\tt q} \in \mathbb{N}^m : }\atop{|{\tt q}| \in \{k,k-2,\ldots\}}}
  T_{k,{\tt q}} \left( {{\bg{\lambda}}} \right)
  \prod_j {H}_{q_j} \left( x_j \right)
 \end{array}
\]
where, using the same notational conventions for the coefficients as for the Hermite polynomials 
themselves:
\[
 \begin{array}{l}
  {T\!e}_{k,{\tt q}} \left( {{\bg{\lambda}}} \right)
  =
  \left.
  \frac{{k!}}{2^i{\tt q}!i!}
  {\bg{\lambda}}^{{\tt q}}
  \left( \left\| {\bg{\lambda}} \right\|_2^2 - 1 \right)^i
  \right|_{i=\frac{k-|{\tt q}|}{2}} \\

   T_{k,{\tt q}} \left( {{\bg{\lambda}}} \right)
  =
  \left.
  \frac{{k!}}{{\tt q}!i!}
  {\bg{\lambda}}^{{\tt q}}
  \left( \left\| {\bg{\lambda}} \right\|_2^2 - 1 \right)^i
  \right|_{i=\frac{k-|{\tt q}|}{2}}
\end{array}
\]
as previously claimed.

\bibliography{universal}
\bibliographystyle{abbrv}

\end{document}